\documentclass[11pt]{amsart}

\usepackage{graphicx}

\usepackage{amssymb}

\begin{document}

\title[Whitney's Theorem for oscillating
functions] { Whitney's Theorem for oscillating on $\mathbb{R}$
functions
}

\author {Yu. Kryakin}

\address[Yuri Kryakin]{
Institute of Mathematics, University of Wroc{\l}aw,
pl. Grunwaldzki 2/4,
50-384  Wroc{\l}aw,
Poland }

\email{kryakin@list.ru \ kryakin@math.uni.wroc.pl}
\thanks{
The  author  supported by Polish founds
for science in years 2006--2008 (research project )
}


\date{}

\maketitle

\begin{abstract}
We find the order of Whitntey's constants for oscillating
functions.
\end{abstract}




\vskip 0.5 cm
\begin{center}
{\bf 1 Introduction }
\end{center}
\vskip .2cm

\vskip .2cm

Whitney's Theorem  \cite{W57} is the difference analogue of the
following Taylor's estimate:

\vskip .2cm

{\it If $ I:=[0,1]$ and

$$
f^{(k)}(0)=0, \quad k=0,\dots, m-1, \eqno (CTM)
$$

then

$$
\sup_{x \in I} | f (x) | \le (m!)^{-1} \sup_{x \in I} |
f^{(m)}(x)|.
\eqno (TM)
$$
}
\vskip .2cm

{\bf Theorem W} \cite{W57, GKS02}. {\it If
$$
f(k/(m-1))=0, \quad k=0, \dots, m-1, \eqno (CIW)
$$

then
$$
\sup_{x \in I} | f(x) | \le 3 \, \sup_{x, x+mh \in I} | \Delta_h^m f(x)|,
\eqno (IW)
$$

where
$$
\Delta_h^m f(x):=\sum_{j=0}^m (-1)^j \binom mj f(x + (j-m/2)h).
$$
}

The inequality (IW) is a consequence of the following Whitney type
inequality \cite{GKS02}:

\vskip .3cm

{\it If
$$
\int_0^{k/m} f(t) \, dt=0, \quad k=1, \dots, m, \eqno (CMW)
$$

then
$$
\sup_{x \in I}   | f(x) | \le  W_m \sup_{x, x+mh \in I} | \Delta_h^m f(x)
|,
\eqno (MW)
$$
$$
W_m \le 2 + 1/e^2.
$$
}

{\bf The Main Conjecture} on Whitney's constants states that $W_m=1$.  This was proved
for $m \le 8$ \cite{Z02, Z04}.

\vskip .3cm

Historical remarks on this theme one can find in  \cite{K02}.
Note,  that the First Sendov's Conjecture about the constant for
the best approximation by algebraic polynomials \cite{S82}
follows from the Main Conjecture.
The Second Sendov's Conjecture \cite{S82} states that we can put
in (IW) 2 instead of 3.
At the time when Sendov's conjectures have been proposed the estimates
of constants were quite rough. Three in the inequality (IW) is
the result of serious efforts and the proof of (IW) is not simple
\cite{GKS02}.

\vskip 0.5cm

In this note we prove the following Whitney type estimate:

\vskip .2cm
{\it Suppose that  $f$ is a locally integrable on $\mathbb{R}$ .
If
$$
\int_{jh}^{(j+1)h} f(u) \, du =0, \qquad j \in \mathbb{Z}, \qquad h >0,
$$

then
$$
\| f \|:=\sup_{x \in \mathbb{R}} |f(x)|  \le W^*_{2k} \omega_{2k} (f,h),
$$

where
$$
\omega_{2k} (f,h):=\sup_{|t|\le h} \| \Delta_t^{2k} f(\cdot)\|.
$$
}
\vskip .3cm

We shall prove  (in section 2) that
the order of $W^*_{2k}$ can be defined up to the log--factor. Namely,
$$
\frac 1{\binom{2k}{k}} \le W^*_{2k} \le \frac{1 +
H_k}{\binom{2k}{k}}, \qquad H_k:=\sum_{j=1}^k 1/j.  \eqno (WK)
$$
However,  the question about exact constant $W^*_{2k}$
is, probably, extremely difficult.

\vskip .3cm

Particularly, we do not know the exact constant in the simplest
case $k=1$. The proof of the estimate (see Section 3)
$$
 1/2 + 3/37 < W^*_2 < 1/2 + 1/8 \eqno (W2)
$$
is the evidence of the character of the difficulties that appear.
We have a significant difference from the interval--case here (see (MW)),
 where
the similar approach allows to obtain the sharp results:  $W_m = 1$
for $m \le 8$. \cite{Z02, Z04}
\vskip .3cm

We think that it is important to find the right order of constants
$W^*_{2k}$. This is the question about the factor $1+ H_k \asymp \ln k$.

\vskip .3cm

The main motivation for the study of the order of the constants
$W_{2k}^*$ is a possible connection of Whitney's type inequality
with Jackson--Stechkin  inequality
$$
E_n (f) \le J_k (\delta) \ \omega_k (f, \delta).
$$
In the recent paper \cite{FKS06} the exact order (with respect to  $k$) of
Jackson--Stechkin constants  $J_k(\delta)$ has been found for
$\delta > \pi/n$. The case  $ \delta = \pi/n$ is the most difficult
not only for approximation in $L^\infty$ and $L^p, \ p \ge 1$ \cite{FKS06}.
Even in the case $L^2$, were sharp results were obtained
\cite{C67,V01}, the exact constant for this argument have not
yet been found.

We hope that Whitney's type theorem we suggest may be useful to
determine the order of the Jackson--Stechkin constants in this
principal case.

On the other hand, note that the results of the paper
\cite{FKS06} are  motivated by Whitney's estimate (MW).
Namely, this estimate was essential for creating the main tool of the paper
\cite{FKS06} --- Favard's inequality for differences.

\vskip .5cm

\vskip .5cm
\begin{center}
{\bf 2 Asymptotic estimate }
\end{center}

\vskip .2cm

{\bf Theorem 1.} {\it If  $f$ is a locally integrable on $\mathbb{R}$
and
$$
\int_{jh}^{(j+1)h} f(u) \, du =0,
\qquad j \in \mathbb{Z}, \quad  h>0, \qquad \eqno(CPW)
$$

then
$$
\| f \| \le W^*_{2k} \ \omega_{2k} (f,h),
$$

with
$$
\frac 1{\binom{2k}{k}} \le W^*_{2k} \le \frac{1 +
H_k}{\binom{2k}{k}}, \qquad H_k:=\sum_{j=1}^k 1/j.  \eqno (PW)
$$
}
{\sl Proof}. The main tool for the proof is Steklov's function
$$
 F(x,y):= \frac{1}{y-x} \int_x^y f(u) \, du, \quad F(x,x):=f(x).
$$
This function of two variables has a better
smoothness than the function $f$ in the following sense.
Put

$$
\Delta_{h_1,h_2}^{2k} F(x,y):= \sum_{j=0}^{2k} (-1)^j \binom{2k}{j}
F(x+(j-k)h_1, y+(j-k)h_2).
$$

\noindent
It is easy to check that

$$ \Delta_{h_1,h_2}^{2k}F(x,y) = \int_0^1 \Delta_{h_1 +
t(h_2-h_1)}^{2k} f(x +t(y-x)) \, dt. $$

Note, that we can assume that
$h=1$ and  $\omega_{2k}(f,h)=1$.
The integral representation above implies that
$ | \Delta_{h_1,h_2}^{2k} F(x,y) | \le 1, \quad h_j \in [0,1].$

\vskip .2cm

The idea of the proof is simple. By the conditions  $F(j,k)=0$
we have
$$
f(j)=F(j,j)={\binom{2k}k}^{-1} \Delta_{0,1}^{2k} F(j,j).
$$

\noindent

We get the lower estimate immediately.
It is sufficient to consider the function equal to $0$, except
one point, where the function is equal to 1.

In the next part of the proof we shall show that when we pass from
integer points to arbitrary ones  we  lose not greater then $ \ln k $.

The proof is based on the following simple identity:

$$ \Delta_{1,0}^{2k} F(0,x) = \Delta_1^{2k} (1/x)  \int_0^x f(u) \, du =
\Delta_1^{2k} (1/x) x F(0,x),
$$

\noindent
which we can rewrite in the form

$$
\int_0^x f(u) \, du = \frac{(x-k) \cdots x \cdots (x+k)}{2k!}
\Delta_{1,0}^{2k} F(0,x). \eqno (E)
$$
\vskip .2cm

\noindent
For the proof of (PW)  we can suppose (by symmetry argument)
that $x \in (0,1/2]$.

\noindent
Let us rewrite  $\Delta_{0,1}^{2k} F(x,x)$ in the form:

$$
\Delta_{0,1}^{2k} F(x,x)= (-1)^k \binom{2k}k f(x) + R_k(x),
$$
$$
R_k (x) := \sum_{j=-k, j \ne 0}^k (-1)^{k-j} \binom{2k}{k+j} F(x,x+j).
$$

\noindent
We need to show that

$$
|R_k (x)| \le H_k.
$$

\noindent
Since
$$
F(x,x+j)= \frac 1j \left( \int_x^0 + \int_0^j + \int_j^{j+x}
\right),
$$
and the difference $\Delta_{0,1}^{2k}$ is symmetric, and $\int_0^j =0$,
 we have only the terms with $\int_j^{j+x}$ in  $R_k(x)$:
$$
R_k(x) = \sum_{j=-k, j \ne 0}^k \frac {(-1)^{k-j}} j \binom{2k}{k+j} \int_j^{j+x}.
$$
The identity (E) with the condition $\omega_{2k} (f,1)=1$ give

$$
\left| \int_0^x \right| \le \binom{x+k}{2k} = x \prod_{j=1}^k (1-x^2/j^2)
{\binom{2k}{k}}^{-1} \le \frac 12 {\binom{2k}{k}}^{-1}, \quad x
\in (0, 1/2].
$$

\noindent
Again, by symmetry
$$
\left| \int_j^{j+x} \right| \le \frac 12 {\binom{2k}{k}}^{-1}
$$

\noindent
and the estimate for  $|R_k(x)|$ is proved.

\vskip .2cm

{\sl Remark 1.}  It is clear that for the periodic
functions with oscillation we get the exact values of constants
(since $R_k(x) =0$).

\vskip 0.5 cm

\begin{center}
{\bf 3 The case  $k=1$ }
\end{center}

\vskip .2cm

{\bf Theorem 2.} {\it
$$
1/2 + 3/37 \le W^*_{2} \le 1/2 + 1/8.
$$
}

{\sl Proof.} {\bf 1.} Prove the upper estimate at first.
Note, that in the case $k=1$ we have the better estimate than
for  $k \ge 2$ in the following sense:
$$
|f(j/2)|= \left|{\binom{2}1}^{-1} \Delta_{0,1/2}^{2} F(j/2,j/2)\right| \le 1/2.
$$
As before we can suppose that $x \in (0,1/2)$.

{\bf 1.a.} If  $x \in (0,1/4]$, then

$$
\Delta_{0,1}^2 F(x,x) = - 2 f(x) + \int_1^{1+x} -
\int_{-1}^{-1+x}
$$

\noindent
The inequality  (see  (E))
$$
\left|\int_0^{x}\right| \le \left|\frac {(x-1)x(x+1)} 2 \right|
\le \frac 18 \ \frac {15}{16} < 1/8,
$$
implies
$$
\left| \int_1^{1+x} \right| < 1/8,
\qquad \left| \int_{-1}^{-1+x}\right| < 1/8,
$$
and
$$
|f(x)| \le |f(1/4)| < 1/2 + 1/8.
$$
{\bf 1.b.} If  $x \in (1/4,1/2]$, then
$$
\Delta_{0,1-x}^2 F(x,x) = - 2 f(x) - \frac{1}{1-x} \int_{-1}^{-1+2x}
$$
and
$$
\left| \frac{1}{1-x} \int_{-1}^{-1+2x} \right| \le \left|
\frac{(2x-1)2x(2x+1)}{2(1-x)} \right|_{x=1/4} = 1/4.
$$
\qed

{\sl Remark 2.} \  It is clear that we can prove
a little better estimate by combining
1.a and 1.b. The equation
$$
(1-x)x(1+x) = (2x-1)2x(2x+1)/(1-x)
$$
has the root $x_0 = 2 - \sqrt{3}$ and instead of upper estimate
$1/2 + 1/8 = 0.6250$ we have $0.6244$.

\vskip .2cm
{\bf 2.}  The following example gives the lower estimate.

Consider the function, which is  piecewise continuous,
equals to $0$ on  $(-\infty,
-1)$ and on $ [2, \infty)$ and equals to  $1/2 + 3/37$ at the point $1/4.$
On the intervals  $I_1=[-1, -0.5)$, $I_2=[-0.5, 1)$ (without the point 1/4),
 $I_3=[1, 5/4)$, $I_4=[5/4,3/2)$ and $ I_5=[3/2, 2)$ the function is linear :

$f(-1)=-12/37$, $f(-0.5-)=-6/37$, $f(-0.5)=12/37$, $f(1-)=-6/37$,
$f(1)=12/37$, $f(5/4 -)=15/37$, $f(5/4)=-10/37$, $f(3/2)=-1/37$,
$f(2-)=-7/37$.

\vskip 1cm
\begin{center}
\includegraphics[width=4in,height= 6in, angle=270]{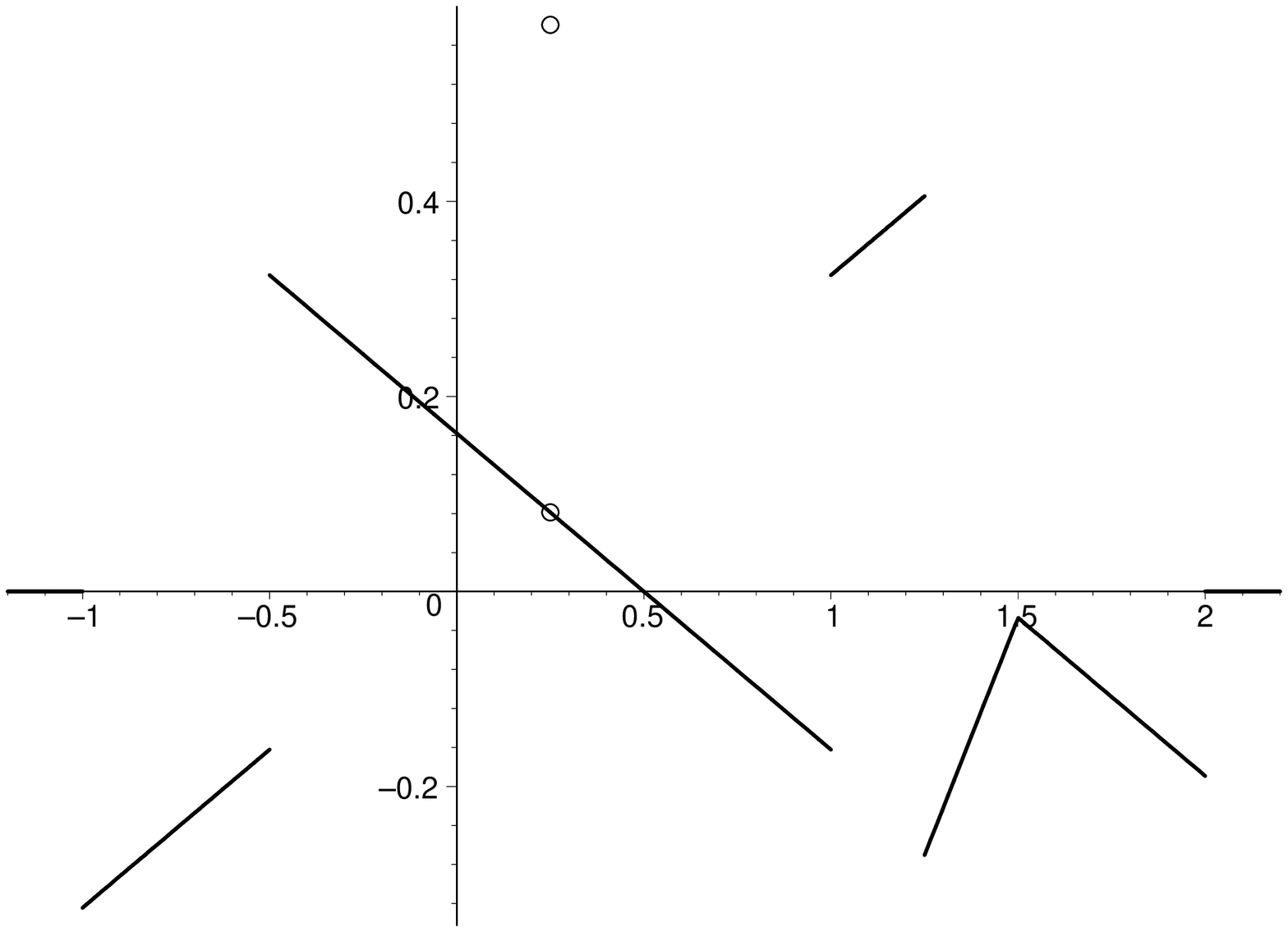}\\
\end{center}
\vskip .2cm

\begin{center} Figure 1. \end{center}

\vskip 0.2cm

The direct computation gives  $\omega_2 (f,1) \le 1$
and  $\int_j^{j+1} = 0, \ j \in \mathbb{Z}$.
Note that we can obtain from this function
by smoothing the example of continuous
function for the lower estimate.
\qed

\end{document}